\documentclass[12pt,a4paper]{article}
\usepackage[cp1251]{inputenc}
\usepackage{graphicx}
\usepackage{amssymb}
\usepackage{amsmath}
\usepackage{amsthm}
\usepackage{euscript}
\numberwithin{equation}{section}
\newtheorem{theorem}{Theorem}
\newtheorem{lemma}{Lemma}

\newtheorem{rem}{Remark}

\newtheorem{ex}{Example}

\newcommand{\RP}[0]{\mathbb{RP}}

\newcommand{\be}[0]{\begin{equation}}
\newcommand{\ee}[0]{\end{equation}}
\newcommand{\bez}[0]{\begin{equation*}}
\newcommand{\eez}[0]{\end{equation*}}
\newcommand{\bl}[0]{\begin{lemma}}
\newcommand{\el}[0]{\end{lemma}}

\newcommand{\ep}[0]{$\hspace{\fill} \square$}
\newcommand{\paragraf}[1]{\par
\bigskip{\centerline{\bf #1}}\medskip}

\newcommand{\abs}[1]{\begin{quotation} {\small
 \centerline{{\bf Abstract}}\smallskip
#1}
\end{quotation}}
\author{I.~Shnurnikov\footnote{NRU HSE}}
\title{About the number of connected components in arrangements of hyperplanes.}
\date{}
\begin{document}
\maketitle
\abs{
We consider arrangements of $n$ hyperplanes of codimension one in real projective $d$--dimensional space $\RP^d$. Let us denote by $f_{max}$ the maximal possible number $f$ of connected components of the complement in space $\RP^d$ to an arrangement of $n$ hyperplanes. We prove that for large enough $n$ and every $d\geq 4$ all integers between $2n\sqrt{n}$ and $f_{max}$ could be realized as numbers $f$ for some arrangements of $n$ hyperplanes in the space $\RP^d$. This fact was known before for $d=2,3$.}

\paragraf{Introduction}
%The main result of the present paper is the following

For arrangements  of $n$ hyperplanes of codimension one in real projective $d$--dimensional space $\RP^d$ Zaslavsky counted the number of connected components of the complement to the union of hyperplanes in terms of lattice of intersections. 
Orlic and Solomon \cite{Orlic} proved that the region number in hyperplane arrangements equals to the cohomology ring dimension of the complement to complexified arrangement. N.~Martinov \cite{Martinov} founded the sets of region numbers in real projective plane arrangements of lines and arrangements of pseudolines.

Let us denote by $C_n^k$ the binomial coefficient $\frac{n!}{k!(n-k)!}$.
Let us denote by $f_{max}=f_{max}(n,d)$ the maximal possible number $f$ of connected components of the complement in space $\RP^d$ to an arrangement of $n$ hyperplanes. It is known that
$$
f_{max}=\sum_{i=0}^{d}C_{n-1}^i
$$
if $n\geq d+1$.
We will use below denotations $n \rightarrow \infty, d\geq 4$ for the number of hyperplanes and the dimension of projective space.

\bigskip 
\paragraf{Main part.}
\begin{lemma}
% count 1
Let we have an arrangement of $n$ hyperplanes in $RP^d$ in general position. We move $k$ hyperplanes so that they will have a common point and this is the only difference between new arrangement and the arrangement in general position. Then the number of regions will decrease by $C_{k-1}^d$.
\end{lemma}

\begin{lemma}
%2
Let $f=f_{max}-t,$ where integer $0 \leq t \leq C_u^2$ and $2d^2<u< \sqrt{n}$. Then there is an arrangement of $n$ hyperplanes in $RP^d$ with $f$ regions, such that it differs from the arrangement in general position in at most $u^2$ hyperplanes and the only difference is some intersection points, which are incident to more then $d$ hyperplanes. 
\end{lemma}
\proof Due to lemma 1 it is enough to notice that every $0 \leq t \leq C_u^2$ could be represented as $\sum_{i=d}^{u-1}n_iC_{i}^d$,
where integers $n_i\leq \frac{i+1}{i+1-d}$ for $i \leq 2d-1$ and $n_i=1$ for $2d \leq i \leq u-1$. It is easy to count that we need to move at most $2d(n_d+\dots+n_2d)+C_u^2<u^2$ hyperplanes.
\ep 	

\begin{lemma}
%3
If $m\rightarrow \infty$ then $C_m^{d-1} \leq C_{m'}^d$ where 
$m'= m^{\tfrac{d-1}{d}}(1+o(1))$.
\end{lemma}

\begin{theorem}
Let integer $0 \leq N \leq C_{n-an^{b}}^d$ where $b=\tfrac{d-1}{d}$ and $a=(d+1)(1+o(1))$ for $n \rightarrow \infty$. 
Then for every $f=f_{max}-N$ there is an arrangement of $n$ hyperplanes in $RP^d$ with $f$ regions, such that all its intersections of dimension at least one are the same as in the arrangement in general position and the only difference is some intersection points, which are incident to more then $d$ hyperplanes. 
\end{theorem}

\proof Let us represent $N=C_m^d+t,$ where $0 \leq t\leq C_{m}^{d-1}$ and $m \leq n-an^{b}$. Then we will use $m+1$ hyperplanes by lemma 1 to reduce the number $f$ from $f_{max}$ to $f_{max}-C_m^d$. Now we will do the same for $N'=t$.
By lemma 3 on the next step we need $m'+1=m^{\tfrac{d-1}{d}}(1+o(1))$ hyperplanes. We will do it until we get $m^{(d)}=m^{\left(\tfrac{d-1}{d}\right)^d}(1+o(1))$.
At the end we use lemma 2. It is easy to count that we need at most 
$$
m+m^{\left(\tfrac{d-1}{d}\right)}+m^{\left(\tfrac{d-1}{d}\right)^2}+\dots+
^{\left(\tfrac{d-1}{d}\right)^d}+ m^{\left(\tfrac{d-1}{d}\right)^{2d}}<n
$$
hyperplanes.
\ep 

For arrangengement $A$ by $A^{(k)}$ we will denote the arrangement, which is produced by hyperplanes of $A$ on the plane in general position of codimension $k$. So $A^{(0)}=A$ and $A^{(d-1)}$ is a set of $n$ points.

\begin{theorem}
For every integer $f_{max(n,d-1)}\leq f\leq f_{max}(n,d)$ there is an arrangement of $n$ hyperplanes in $\RP^d$ with $f$ regions.
\end{theorem}

\proof Let us consider arrangements $A'$ of $n-k$ hyperplanes in $\RP^{d-1}$ which appeared in theorem 1. Let us take a conus $A $ on the arrangement $A$. Let us consider $k$ new hyperplanes which are in general position to $A$ and to each other. Then the total arrangements of $n$ hyperplanes gives $(k+1)f(A)+\sum_{i=2}^{d-1}C_k^i f(A^{(i)})+C_k^d$. As arrangements $A$ were constructed by theorem 1 then $f(A^{(i)})$ is the same for $i>1$. Then we could move $k$ last hyperplanes (changing generality of position only in intersection points) to reduce $f$ by $0 \leq q \leq k.$ Here we need $k \geq O(d^2)$. Now we take $O(d^2) \leq k \leq n^{1-\frac 1d}$ and obtain any integer 
$$
f_{max}(n,d-1)\leq f\leq f_{max}(n,d)-C_{n-an^{b}}^d
$$ (where 
$b=\tfrac{d-1}{d}$ and $a=(d+1)(1+o(1))$ for $n \rightarrow \infty$)
as a number of regions in an arrangement of $n$ hyperplanes in $\RP^{d}$. Theorem 1 completes the proof.
\ep 

\begin{theorem}
For every integer $\sqrt{2n}n \leq f\leq f_{max}(n,d)$ there is an arrangement of $n$ hyperplanes in $\RP^d$ with $f$ regions.
\end{theorem}
\proof We use induction on $n$ and theorems 1 and 2. \ep

\end{document}